\documentclass[11pt,a4paper]{article}
\usepackage{graphicx}
\usepackage{amsmath}
\usepackage{amssymb}
\usepackage{amsthm}
\usepackage{cite}
\usepackage{subfigure}

\providecommand{\keywords}[1]
{
	\small	
	\textbf{\textit{Keywords---}} #1
}

\newcommand{\bm}{\boldsymbol}
\newcommand{\vect}[1]{\mathbf{#1}}
\newcommand{\mat}[1]{\mathbf{#1}}

\newcommand{\E}[1]{{\mathbb{E}}\{#1\}}

\newcommand{\etr}[1]{{\mathrm{etr}}\left\{#1\right\}}

\renewcommand{\det}[1]{|#1|}
\newcommand{\Prob}[1]{{\mathbb{P}}\left(#1\right)}

\newcommand{\alphaML}{\alpha_{\text{\tiny{ML}}}}
\newcommand{\loss}{\rho}
\newcommand{\lossGaussian}{\loss_{\text{\tiny{Gaussian}}}}
\newcommand{\lossStudent}{\loss_{\text{\tiny{Student}}}}
\newcommand{\betaGaussian}{\beta_{\text{\tiny{Gaussian}}}}

\newcommand{\ttilde}{\tilde{t}}
\newcommand{\ttildeGaussian}{\ttilde_{\text{\tiny{Gaussian}}}}

\newcommand{\xtilde}{\tilde{x}}
\newcommand{\Pfa}{P_{fa}}

\newcommand{\elast}{\vect{e}_{N}}
\newcommand{\vn}{\vect{n}}

\newcommand{\vs}{\vect{s}}
\newcommand{\vt}{\vect{t}}
\newcommand{\vv}{\vect{v}}
\newcommand{\vw}{\vect{w}}
\newcommand{\vwamf}{\vw_{\text{\tiny{amf}}}}
\newcommand{\vwopt}{\vw_{\text{\tiny{opt}}}}
\newcommand{\vx}{\vect{x}}
\newcommand{\vxtilde}{\tilde{\vx}}
\newcommand{\vzero}{\vect{0}}

\newcommand{\F}{\mat{F}}

\newcommand{\I}{\mat{I}}
\newcommand{\eye}[1]{\I_{#1}}
\newcommand{\N}{\mat{N}}
\newcommand{\Q}{\mat{Q}}

\renewcommand{\S}{\mat{S}}
\newcommand{\T}{\mat{T}}

\newcommand{\W}{\mat{W}}
\newcommand{\X}{\mat{X}}
\newcommand{\Xbar}{\bar{\X}}
\newcommand{\Y}{\mat{Y}}

\newcommand{\Mzero}{\mat{0}}

\newcommand{\mDelta}{\bm{\Delta}}

\newcommand{\mSigma}{\bm{\Sigma}}
\newcommand{\mSigmat}{\bm{\Sigma}_{t}}
\newcommand{\mOmega}{\bm{\Omega}}

\newcommand{\CN}[2]{\tilde{\mathrm{N}}\left(#1,#2\right)}
\newcommand{\vCN}[3]{\tilde{\mathrm{N}}_{#1}\left(#2,#3\right)}
\newcommand{\mCN}[4]{\tilde{\mathrm{N}}_{#1}\left(#2,#3,#4\right)}
\newcommand{\vCT}[4]{\tilde{\mathrm{T}}_{#1}\left(#2,#3,#4\right)}
\newcommand{\mCT}[5]{\tilde{\mathrm{T}}_{#1}\left(#2,#3,#4,#5\right)}
\newcommand{\CW}[3]{\tilde{\mathrm{W}}_{#1}\left(#2,#3\right)}

\newcommand{\Cchisquare}[2]{\tilde{\chi}^{2}_{#1}(#2)}
\newcommand{\CF}[3]{\tilde{\mathrm{F}}_{#1}(#2,#3)}

\newcommand{\dist}{\overset{d}{=}}

\newcommand{\onehalf}{\frac{1}{2}}
\newcommand{\SNR}{\mathrm{SNR}}

\newcommand{\Beta}[2]{B_{#1,#2}}

\begin{document}
\title{On the distributions of some statistics related to adaptive filters trained with $t$-distributed samples}
\author{Olivier Besson\thanks{The author is with Universit\'{e} de Toulouse, ISAE-SUPAERO, 10 avenue Edouard Belin, 31055 Toulouse, France. Email: olivier.besson@isae-supaero.fr}}
\date{March 2021}
\maketitle
\begin{abstract}
In this paper we analyse the behaviour of adaptive filters or detectors when they are trained with $t$-distributed samples rather than Gaussian distributed samples. More precisely we investigate the impact on the distribution of some relevant statistics including the signal to noise ratio loss and the Gaussian generalized likelihood ratio test. Some properties of partitioned complex $F$ distributed matrices are derived which enable to obtain statistical representations in terms of independent chi-square distributed random variables. These representations are compared with their Gaussian counterparts and numerical simulations illustrate and quantify the induced degradation. 
\end{abstract}

\section{Introduction}
Estimating the amplitude $\alpha$ or detecting the presence of a known signal $\vv$ from a noise corrupted version $\vx=\alpha\vv+\vn$ is a recurrent problem in numerous applications including radar where $\vv$ stands for the space and/or time signature of a potential target and $\vn$ gathers disturbance sources, mostly clutter and thermal noise \cite{Ward94,Melvin12,Richards14}. When the noise $\vn$ follows a complex matrix-variate Gaussian distribution with zero mean and covariance matrix $\mSigma$ the maximum likelihood estimate (MLE) of $\alpha$ writes $\alphaML = \vwopt^{H} \vx$ with $\vwopt= (\vv^{H} \mSigma^{-1} \vv)^{-1} \mSigma^{-1} \vv$. This optimal filter $\vwopt$ is also obtained as the solution to the following minimization problem
\begin{equation}
	\underset{\vw}{\min} \vw^{H}\mSigma\vw \text{ subject to } \vw^{H}\vs=1	
\end{equation}
In other words this filter minimizes the output power under the constraint that the signal of interest goes unscathed through the filter.  Note that this interpretation holds irrespective of the distribution of $\vn$. Since $\mSigma$ is usually unknown a set of training samples is used which, in the best case, share the same distribution as $\vn$. In the Gaussian framework $\mSigma$ is substituted for the sample covariance matrix (SCM) $\S=\X\X^{H}$ where $\X$ is the training samples data matrix, on the rationale that $\S$ is up to a scaling factor the MLE of $\mSigma$. Proceeding this way results in $\vwamf = (\vv^{H} \S^{-1} \vv)^{-1} \S^{-1} \vv$ which is usually referred to as the adaptive matched filter \cite{Robey92}. For any filter $\vw$ a classical figure of merit is the signal to noise ratio (SNR) loss which is defined as 
\begin{equation}\label{def_SNRloss}
	\loss(\vw) = \frac{\SNR(\vw)}{\SNR(\vwopt)}	 = \frac{|\vw^{H}\vv|^{2}}{(\vv^{H}\mSigma^{-1}\vv)(\vw^{H}\mSigma\vw)}
\end{equation}
and corresponds to the ratio of the SNR obtained with $\vw$ to that obtained with $\vwopt$.  In the sequel, we concentrate on the SNR loss of $\vwamf$, which we will denote as $\loss$ and is given by
\begin{equation}\label{SNRloss_AMF}
	\loss = \loss(\vwamf) = \frac{(\vv^H \S^{-1} \vv)^{2}}{(\vv^{H}\mSigma^{-1}\vv)(\vv^{H}\S^{-1}\mSigma\S^{-1}\vv)}
\end{equation}
Assuming a Gaussian distribution for $\X$, it has been shown that $\rho$ is beta distributed with parameters that depend only on the size of the observations $N$ and the number of training samples $K$ \cite{Reed74,Khatri87}. A similar beta distribution with different parameters is obtained when persymmetry is exploited \cite{Liu16}. 

Unfortunately for some applications it may not be possible to dispose of Gaussian distributed training samples as the latter have possibly a heavier distribution tail. This is often the case in radar applications where the main source of noise is the clutter and the latter is generally non Gaussian \cite{Farina97,Billingsley99,Conte05}. Therefore, it becomes of interest to study what happens when training samples are no longer Gaussian distributed. This is the aim of this paper where we assume that $\X$ follows a matrix-variate complex $t$ (Student) distribution and we study the impact on the distribution of some random variables commonly used in adaptive filtering and detection, including the SNR loss. As we shall see later, the matrix-variate complex Student distribution also appears naturally when training samples exhibit a particular case of covariance mismatch. In this paper we derive stochastic representations of relevant statistics in terms of independent random variables following a complex chi-square distribution. These representations rely on some properties of partitioned complex $F$ distributed matrices. They allow quick insights into the impact of mismatched training samples.

We note that in the literature the impact of mismatch on adaptive filters or adaptive detectors has been extensively studied, with two main types of mismatch considered. The first concerns a mismatch on the SoI signature $\vv$, see e.g., \cite{Boroson80,Kelly89b,Kalson92,Bose96,Bandiera09,Liu19}. Alternatively, researchers have studied the case where the covariance matrix of $\X$ differs from that of the data to be filtered or the data under test \cite{Richmond00b,Raghavan19,Besson20e,Besson21}. A possible combination of the two mismatches is addressed in \cite{Blum00,McDonald00}. The situation considered herein is different as the mismatch concerns the training samples distribution.

Before proceeding we state the notations used in this paper concerning matrix-variate distributions (MVD). References \cite{Muirhead82,Gupta00} provide a very comprehensive overview of real-valued MVD. For their extension to complex-valued MVD  we refer to e.g. \cite{Goodman63,Khatri65,Krishnaiah76,Mathai05} where most of the distributions considered below are studied. In the sequel we note $\mCN{p,n}{\Xbar}{\mSigma}{\mOmega}$ the complex matrix-variate distribution whose probability density function (p.d.f.) is $p(\X) = \pi^{-pn} \det{\mSigma}^{-n} \det{\mOmega}^{-p} \etr{-\mSigma^{-1}(\X-\Xbar)\mOmega^{-1}(\X-\Xbar)^{H}}$. When $\Xbar=\E{\X}=\Mzero$ the matrix $\S=\X\X^{H} \dist \CW{p}{n}{\mSigma}$ follows a complex Wishart distribution with p.d.f. $p(\S) \propto \det{\mSigma}^{-n} \det{\S}^{n-p} \etr{-\mSigma^{-1}\S}$ where $\propto$ means ``proportional to''. The complex matrix-variate $t$ distribution is denoted by $\mCT{p,n}{\nu}{\Xbar}{\mSigma}{\mOmega}$ and its p.d.f is given by $p(\X) \propto \det{\mSigma}^{-n} \det{\mOmega}^{-p} \det{\I_{p}+\mSigma^{-1}(\X-\Xbar)\mOmega^{-1}(\X-\Xbar)^{H}}^{-(\nu+p+n-1)}$. It is the distribution of $\X = \Xbar + (\W^{-1/2})^{H}\Y$ where $\W \dist \CW{p}{\nu+p-1}{\mSigma^{-1}}$ is independent of $\Y \dist \mCN{p,n}{\mat{0}}{\I_{p}}{\mOmega}$. $\W^{1/2}$ denotes any square-root of $\W$ while $\W^{\frac{1}{2}}$ will stand for its unique Hermitian square-root. If $\S_{i} \dist \CW{p}{n_{i}}{\mSigma}$, $i=1,2$, then $\F=\S_{1}^{\onehalf} \S_{2}^{-1}\S_{1}^{\onehalf}$ follows a complex matrix-variate $F$ distribution with p.d.f. $p(\F) \propto \det{\F}^{n_{1}-p} \det{\I_{p}+\F}^{-(n_{1}+n_{2})}$ and we note $\F \dist \CF{p}{n_{1}}{n_{2}}$. The complex chi-square distribution with $q$ degrees of freedom and non-centrality parameter $\delta$ will be denoted as $\Cchisquare{q}{\delta}$.

\section{Analysis of SNR loss with Student distributed training samples}
In the sequel we assume that $K$ training samples are available and distributed according to $\X \dist \mCT{N,K}{\nu-N+1}{\Mzero}{\mu\mSigma}{\I_{K}}$ so that their p.d.f is given by
\begin{equation}\label{pdf_Student}
	p(\X) \propto \det{\mu\mSigma}^{-K} \det{\I_{N}+(\mu\mSigma)^{-1}\X\X^{H}}^{-(\nu+K)}	
\end{equation}
As explained above, there are situations where the training samples are not Gaussian distributed, e.g., in radar applications where the dominant part of the noise, namely the clutter, is often non Gaussian and well modelled by the class of compound-Gaussian distributions, of which the Student distribution is a member.  Other applications have to deal with non Gaussian data and therefore it is of interest to investigate what happens when a filter is trained with samples that no longer follow a Gaussian distribution but rather a Student distribution. Note that the SNR loss, as given in \eqref{def_SNRloss}-\eqref{SNRloss_AMF}, does not depend on the distribution of the data to be filtered, it just requires that their covariance matrix is $\mSigma$. A second motivation for the use of the Student distribution is the following. Assume that the training samples are Gaussian distributed but have a covariance matrix $\mSigmat$ that is different from $\mSigma$, say with no loss of generality $\mSigmat=\mSigma^{1/2} \W^{-1} (\mSigma^{1/2})^{H}$ for some positive definite matrix $\W$. This is the case for instance in non homogeneous environments in radar applications. We can thus assume that $\X | \W \dist \mCN{N,K}{\mat{0}}{\mSigma^{1/2} \W^{-1} (\mSigma^{1/2})^{H}}{\I_{K}}$. In \cite{Besson20e} we analysed the distribution of the SNR loss for fixed and arbitrary $\W$. We showed that it can be written as a quadratic form in normal or Student random variables and we proposed approximations of them. Now the matrix $\W$ may be considered as a random matrix and, if we assume a conjugate prior $\W \dist \CW{N}{\nu}{\mu^{-1}\I_{N}}$, then the marginal distribution of $\X$ is given by \eqref{pdf_Student}. In other words, the statistical model used herein results from a Bayesian model of  covariance mismatch where the samples used to train the filter do not share the same covariance matrix as the samples to be filtered.  Therefore the model used in this paper covers the two cases described above. Note that the smaller $\nu$ the more heavy-tailed is the Student distribution. 

The sample covariance matrix $\S=\X\X^{H}$ is still, up to a scaling factor, the MLE of $\E{\X\X^{H}}=K(\nu-N)^{-1}\mu \mSigma$ and thus can still be used to design the adaptive filter $\vwamf = (\vv^{H} \S^{-1} \vv)^{-1} \S^{-1} \vv$ whose SNR loss we are interested in. First let us note that $\X \dist (\mu\mSigma)^{\onehalf} \W_{\nu}^{-\onehalf}\N$ where $\W_{\nu} \dist \CW{N}{\nu}{\I_{N}}$ is independent of $\N \dist \mCN{N,K}{\mat{0}}{\I_{N}}{\I_{K}}$ \cite{Gupta00} so that
\begin{align}
	\S 	\dist \mu \mSigma^{\onehalf} \W_{\nu}^{-\onehalf} \W_{K} \W_{\nu}^{-\onehalf} \mSigma^{\onehalf} =  \mu \mSigma^{\onehalf} \F^{-1} \mSigma^{\onehalf}
\end{align}
where $\W_{K} \dist \CW{N}{K}{\I_{N}}$. It follows that $\F=\W_{\nu}^{\onehalf} \W_{K}^{-1} \W_{\nu}^{\onehalf} \dist \CF{N}{\nu}{K}$ \cite{Khatri65,Olkin64}. Therefore the SNR loss can be represented as
\begin{align}
	\loss &= \frac{\SNR(\vwamf)}{\SNR(\vwopt)} = \frac{(\vv^H \S^{-1} \vv)^{2}}{(\vv^{H}\mSigma^{-1}\vv)(\vv^{H}\S^{-1}\mSigma\S^{-1}\vv)} \nonumber \\
	&\dist  \frac{(\vv^H \mSigma^{-\onehalf} \F \mSigma^{-\onehalf}  \vv)^{2}}{(\vv^{H}\mSigma^{-1}\vv)(\vv^{H}\mSigma^{-\onehalf}  \F^{2} \mSigma^{-\onehalf}  \vv)} \nonumber \\
	&\dist  \frac{(\vv^H \mSigma^{-\onehalf} \Q \F \Q^{H} \mSigma^{-\onehalf} \vv)^{2}}{(\vv^{H}\mSigma^{-1}\vv)(\vv^{H}\mSigma^{-\onehalf} \Q \F^{2} \Q^{H}\mSigma^{-\onehalf}  \vv)}
\end{align}
for any unitary matrix $\Q$ since $\F$ and $\Q^{H}\F\Q$ have the same distribution. Let us choose $\Q$ such that $\Q^{H} \mSigma^{-\onehalf}  \vv = (\vv^{H}\mSigma^{-1}\vv)^{1/2} \elast$ where $\elast = \begin{bmatrix} 0 & \ldots & 0 & 1 \end{bmatrix}^{T}$. Partitioning $\F$ as
\begin{equation}\label{partition_F}
	\F = \begin{pmatrix} \F_{11} & \F_{12} \\ \F_{21} & F_{22} \end{pmatrix}	
\end{equation}
where $\F_{11}$ is $(N-1) \times (N-1)$, we arrive at
\begin{align}
	\loss &\dist \frac{(\elast^{H}\F\elast)^{2}}{\elast^{H}\F^{2}\elast} = \frac{F_{22}^{2}}{F_{22}^{2}+\F_{21}\F_{12}} =  \frac{1}{1+\vt_{12}^{H}\vt_{12}}
\end{align}
with $\vt_{12} = \F_{12}F_{22}^{-1}$. As shown in \ref{app:F}, one has
\begin{equation}\label{rep_t12}
	\vt_{12} \dist (1+F_{22}^{-1})^{1/2} \frac{\vn_{12}}{\sqrt{\gamma_{12}}}
\end{equation}
where $F_{22}$, $\vn_{12}$ and $\gamma_{12}$ are independent with $\vn_{12} \dist \vCN{N-1}{\vzero}{\I_{N-1}}$ and
\begin{align}\label{rep_F22}
	\gamma_{12} \dist \Cchisquare{K-N+2}{0}; \quad F_{22} \dist \frac{\Cchisquare{\nu}{0}}{\Cchisquare{K-N+1}{0}}
\end{align}
It ensues that the SNR loss admits the following representation
\begin{equation}\label{rep_SNRloss_Student}
	\lossStudent \dist \left[1 + \left(1+\frac{\Cchisquare{K-N+1}{0}}{\Cchisquare{\nu}{0}}\right) \frac{\Cchisquare{N-1}{0}}{\Cchisquare{K-N+2}{0}}\right]^{-1}
\end{equation}
which provides a simple and convenient expression as a function of independent chi-square distributed random variables. This should be compared to its counterpart when $\X$ is Gaussian distributed, namely
\begin{equation}\label{rep_SNRloss_Gaussian}
	\lossGaussian \dist \left[1 + \frac{\Cchisquare{N-1}{0}}{\Cchisquare{K-N+2}{0}}\right]^{-1}
\end{equation}
Clearly the SNR loss is likely to take lower values in the Student case than in the Gaussian case and we recover that the two representations are equivalent as $\nu \rightarrow \infty$. Moreover the average value of the term $\Cchisquare{K-N+1}{0} / \Cchisquare{\nu}{0}$ is $(K-N+1)/(\nu-1)$ and hence the difference is expected to increase as $K$ increases. The representation in \eqref{rep_SNRloss_Student} also allows to derive (see \ref{app:mean_SNRloss_Student}) the SNR loss p.d.f. which is given in equation \eqref{p(rho)} as well as its mean value which writes 
\begin{align}\label{mean_rho_Student}
	&\E{\lossStudent} = \frac{\nu(K-N+2)}{(\nu+K-N+1)(K+1)} \nonumber \\
	&\times {\;}_{3}F_{2}(1,K-N+3,K-N+1;K+2,\nu+K-N+2;1)	
\end{align}
to be compared with $\E{\lossGaussian} = (K-N+2)/(K+1)$. 

\begin{figure}[h]
	\centering
	\subfigure{%
		\includegraphics[width=11cm]{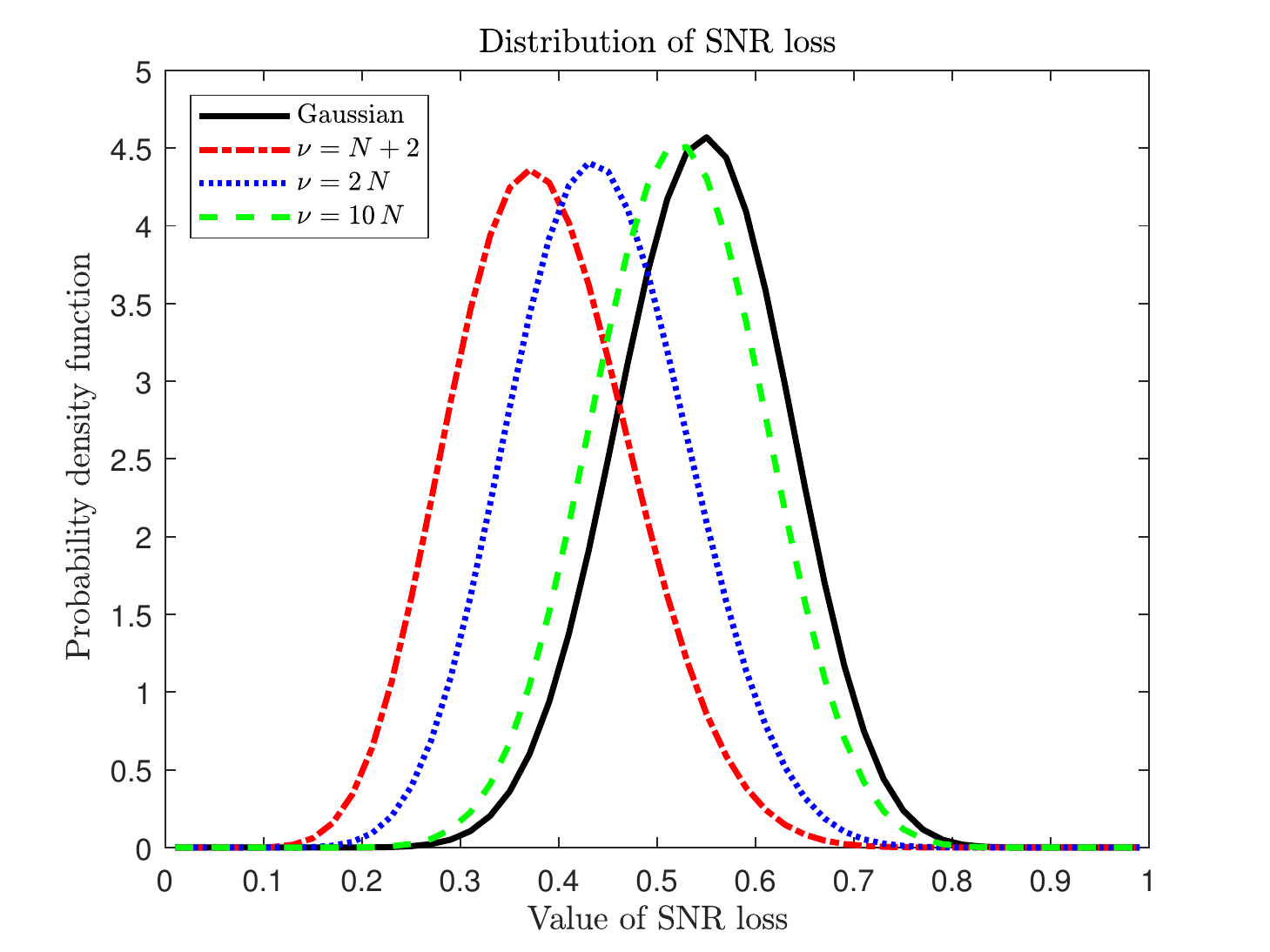}} \\
	\subfigure{%
		\includegraphics[width=11cm]{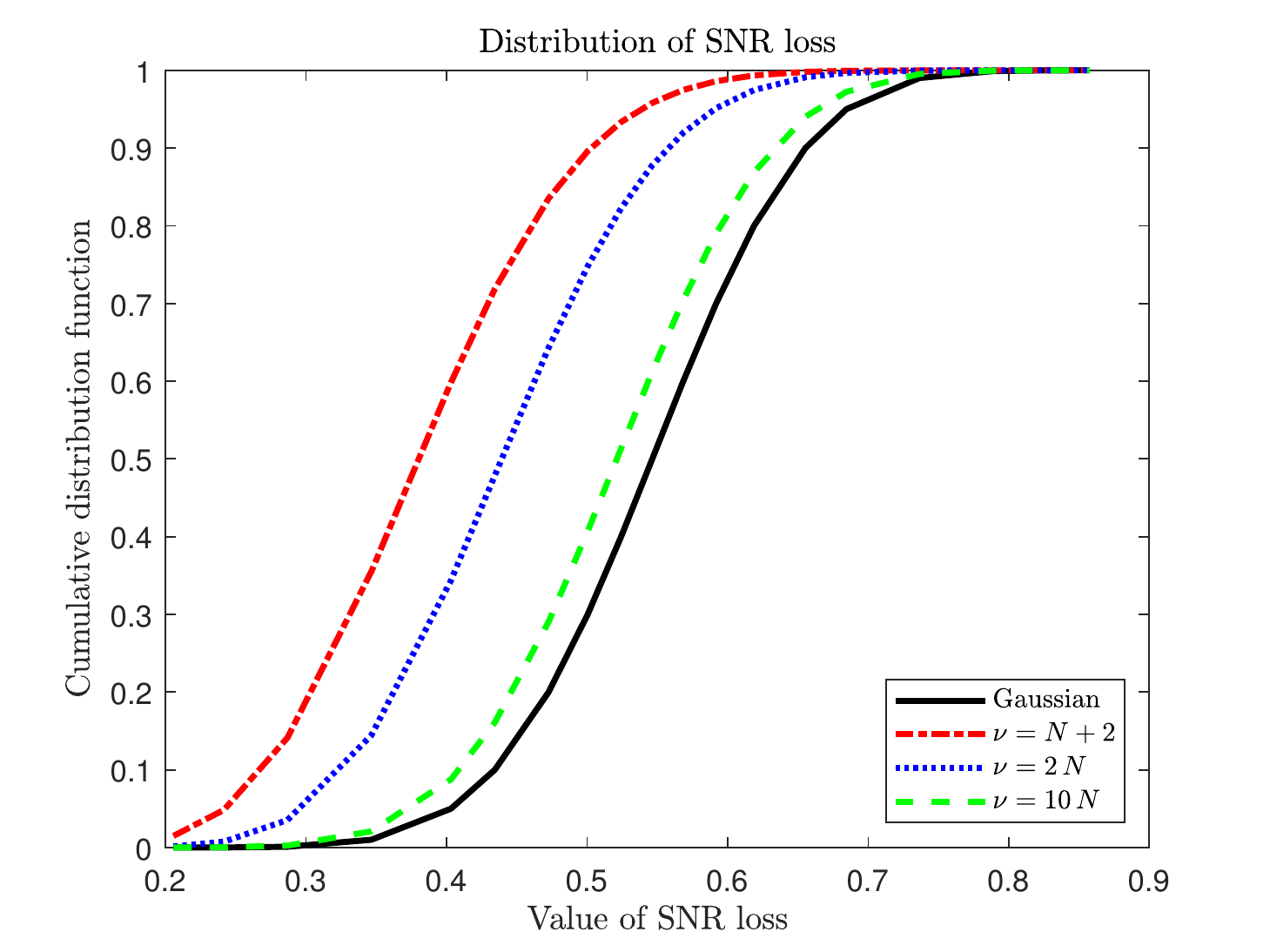}}
	\caption{Probability density function and cumulative distribution function of $\loss$ for various $\nu$. $K=2N$.}
	\label{fig:PDF-CDF_rho_K=32_mean_iW=1}
\end{figure}

\begin{figure}[h]
	\centering
	\subfigure{%
		\includegraphics[width=11cm]{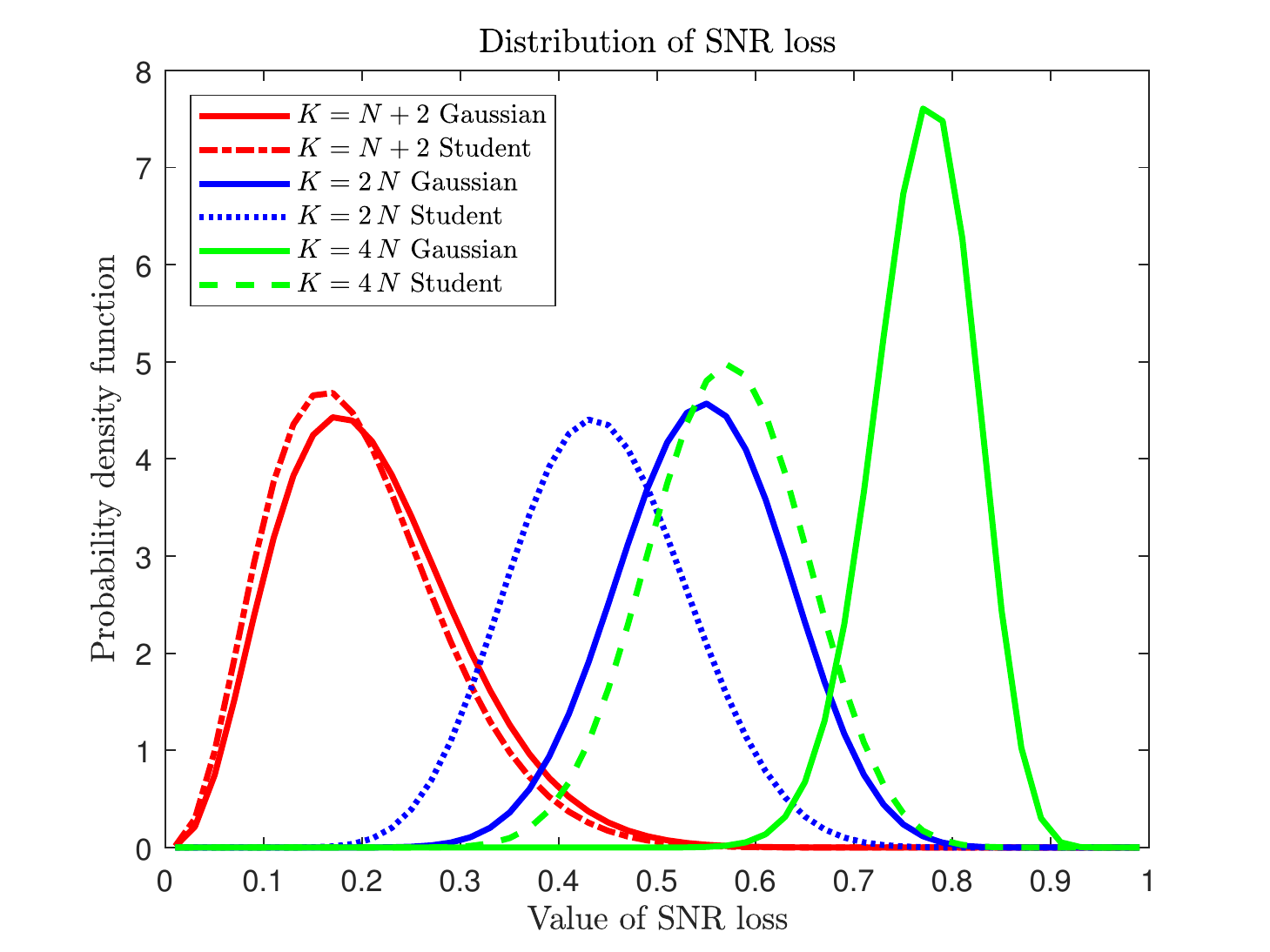}} \\
	\subfigure{%
		\includegraphics[width=11cm]{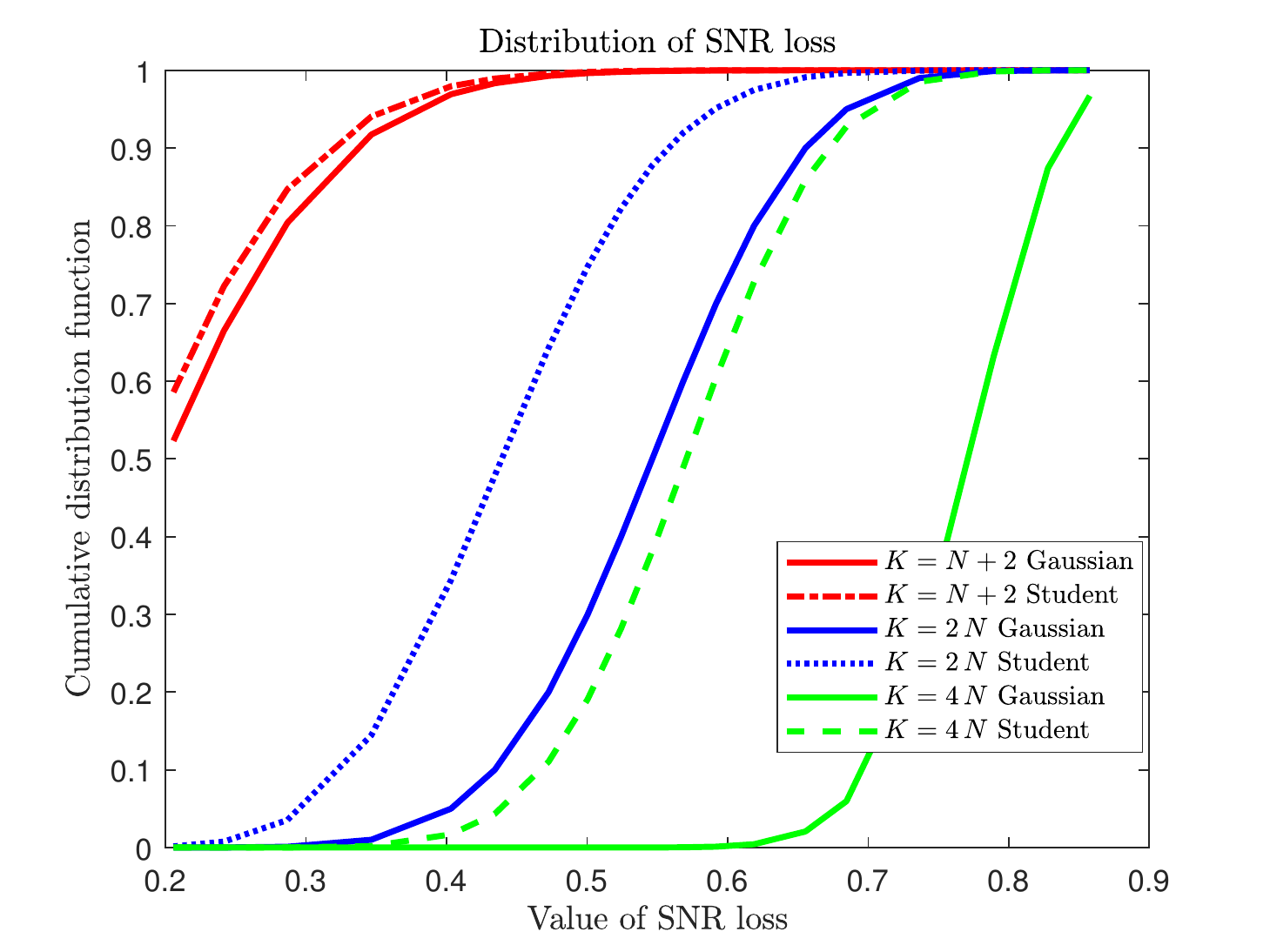}}
	\caption{Probability density function and cumulative distribution function of $\loss$ for various $K$. $\nu=2N$.}
	\label{fig:PDF-CDF_rho_nu=32_mean_iW=1}
\end{figure}

\begin{figure}[h]
	\centering
	\includegraphics[width=11cm]{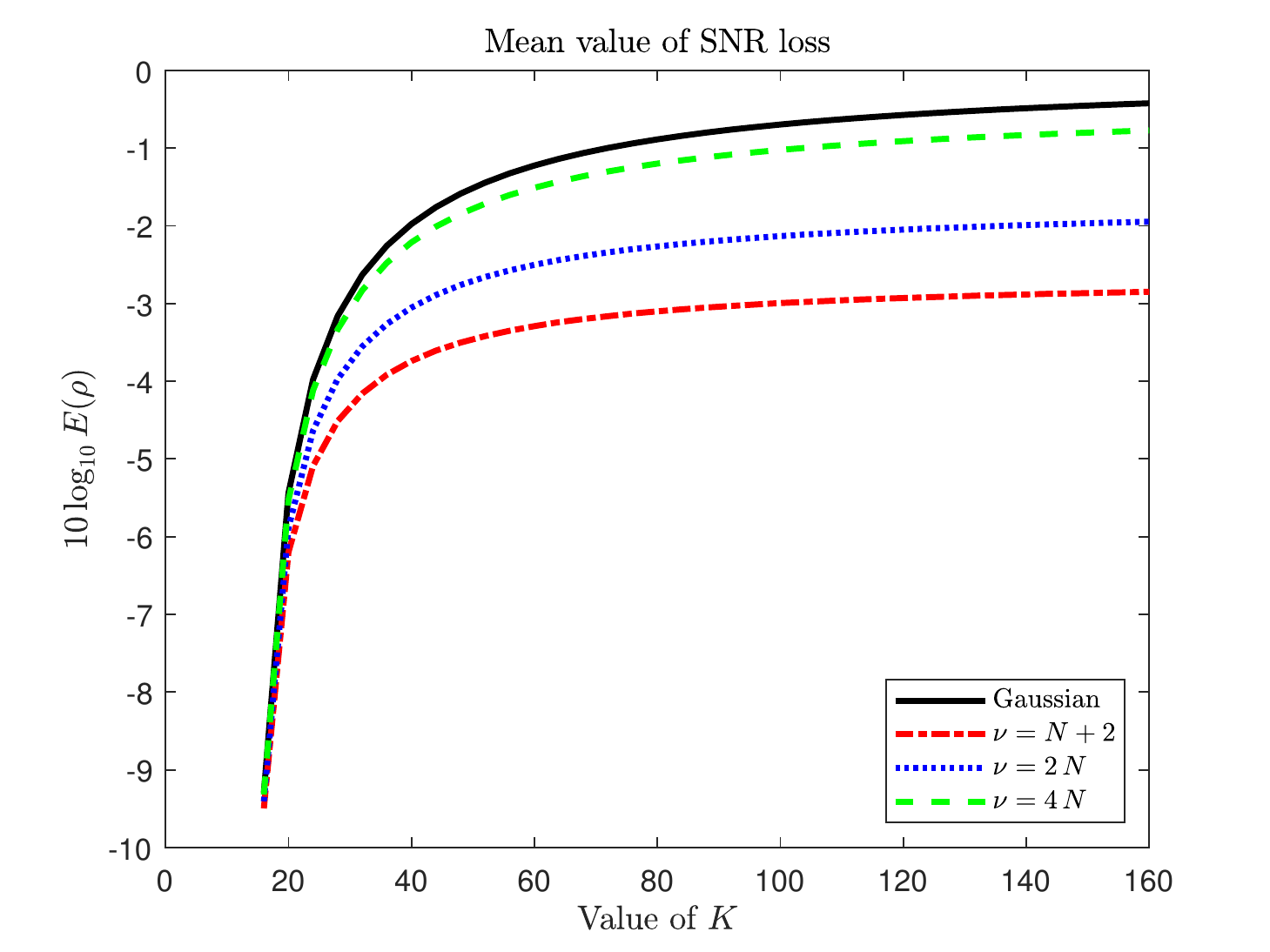}
	\caption{Average value of $\lossStudent$ versus $K$ for various $\nu$.}
	\label{fig:mean_rho_vs_K_mean_iW=1}
\end{figure}

We now provide numerical evaluation of the difference between the distribution of the SNR loss obtained with Gaussian training samples and that obtained with Student training samples. Through preliminary simulations we checked that the distribution of the SNR loss obtained from the representation in \eqref{rep_SNRloss_Student} coincides with the distribution obtained when one generates snapshots from \eqref{pdf_Student}, computes $\vwamf$ and its SNR loss in \eqref{SNRloss_AMF}. We consider a scenario with $N=16$ and $\mu$ is chosen equal to $\nu-N$. We first look at the influence of $\nu$ in Figure \ref{fig:PDF-CDF_rho_K=32_mean_iW=1} where we display the p.d.f and the cumulative distribution function (c.d.f.) of $\loss$  for $K=2N$. As can be seen, the impact is rather significant. For instance while $\Prob{\lossGaussian \leq 0.5}=0.3$ we have $\Prob{\lossStudent \leq 0.5}=0.4$, $0.748$ and $0.896$ for $\nu=10N$, $\nu=2N$ and $\nu=N+2$ respectively. This impact depends however on $K$ as illustrated in Figure \ref{fig:PDF-CDF_rho_nu=32_mean_iW=1}. As could be expected from \eqref{rep_SNRloss_Student}, the difference between the Student and the Gaussian cases increases with $K$. For instance for $K=4N$ the probability of having an SNR loss lower than $0.5$ increases from $\Prob{\lossGaussian \leq 0.5}=3.65\, 10^{-6}$ to $\Prob{\lossStudent \leq 0.5}=0.19$, while for $K=2N$ one goes from $\Prob{\lossGaussian \leq 0.5}=0.3$ to $\Prob{\lossStudent \leq 0.5}=0.748$. This is further illustrated in Figure \ref{fig:mean_rho_vs_K_mean_iW=1} where we display the average value of the SNR loss versus the number of snapshots. As can be seen, the larger $K$ the larger the difference between $\E{\lossGaussian}$ and $\E{\lossStudent}$.

\begin{figure}[h]
	\centering
	\includegraphics[width=11cm]{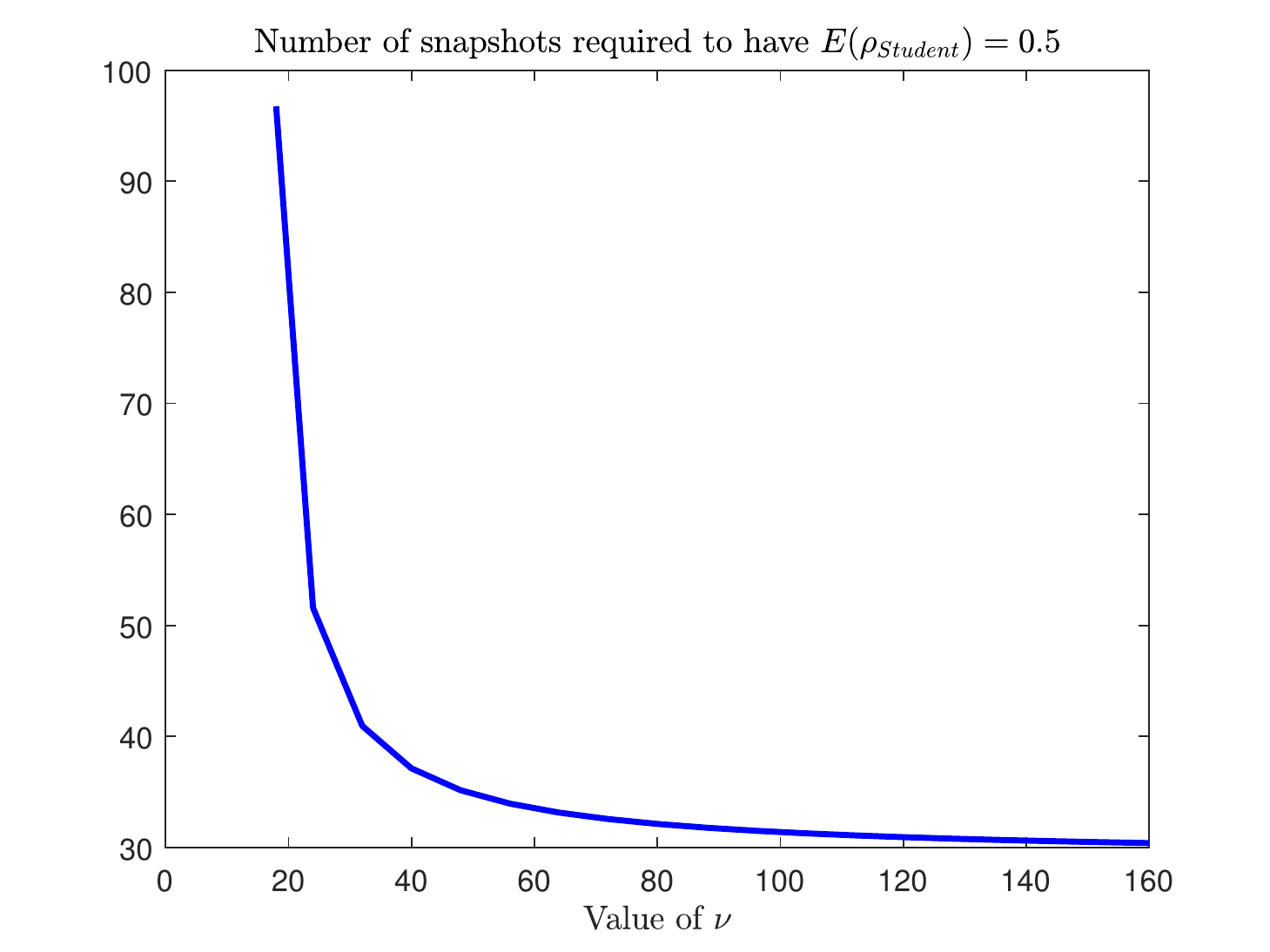}
	\caption{Number of snapshots required to have $\E{\lossStudent}=0.5$ versus $\nu$.}
	\label{fig:find_K_rho_onehalf_Student}
\end{figure}

Another impact concerns the rate of convergence of the adaptive filter which is increased with Student training samples as can be observed in Figure \ref{fig:find_K_rho_onehalf_Student} where we plot the value of $K$ required to have $\E{\lossStudent}=0.5$: clearly the required number of samples decreases when $\nu$ increases, going from $K \simeq 30$ in the Gaussian case to $K \simeq 96$ when $\nu=N+2$.

\section{Distribution of some statistics related to adaptive detection}
We now study the impact of Student distributed training samples for a related problem, namely that of adaptive detection. A very common problem in multichannel processing \cite{Kelly86} is to test $\mathcal{H}_{0}$ versus $\mathcal{H}_{1}$ where
\begin{align}\label{prob_detect}
	\mathcal{H}_{0}: \, &\vx \dist \vCN{N}{\vzero}{\mSigma}; \X \dist \mCN{N,K}{\Mzero}{\mSigma}{\eye{K}} \nonumber \\
	\mathcal{H}_{1}: \, &\vx \dist \vCN{N}{\alpha \vv}{\mSigma}; \X \dist \mCN{N,K}{\Mzero}{\mSigma}{\eye{K}}
\end{align}
The maximal invariant statistic for the detection problem in \eqref{prob_detect} is bi-dimensional \cite{Bose95} and is a one-to-one function of $\beta=1/(1+s_{1}-s_{2})$ and $\ttilde= s_{2}/(1+s_{1}-s_{2})$ where
\begin{equation}
	s_{1}  =\vx^{H}\S^{-1}\vx ; \, s_{2} = \frac{|\vx^{H}\S^{-1}\vv|^{2}}{\vv^{H}\S^{-1}\vv}
\end{equation}
$\beta$ corresponds to the loss factor whose distribution is actually that of $\lossGaussian$ when $\X$ is Gaussian distributed. $\ttilde$ corresponds to Kelly's generalized likelihood ratio test (GLRT) statistic \cite{Kelly86}. Any detector which is a function of $(\beta,\ttilde)$ has a constant false alarm rate property and actually most of the adaptive detectors derived so far can be expressed as a function of $(\beta,\ttilde)$ \cite{Coluccia20}. Therefore, it is of interest to see how the performance of these detectors is affected when the training samples are no longer Gaussian distributed but Student distributed. Note that the impact of a fixed covariance mismatch between $\vx$ and $\X$ with the latter being both Gaussian distributed has been studied in \cite{Richmond00b,Blum00,McDonald00,Raghavan19,Besson21}. In the sequel we consider a distribution mismatch and we assume that $\vx \dist \vCN{N}{\alpha \vv}{\mSigma}$ (where $\alpha$ is possibly equal to zero) and that $\X \dist \mCT{N,K}{\nu-N+1}{\Mzero}{\mu\mSigma}{\I_{K}}$  as before. It is different from assuming that $\begin{bmatrix} \vx & \X \end{bmatrix} \dist \mCT{N,K+1}{\nu-N+1}{\begin{bmatrix} \alpha \vv  & \Mzero \end{bmatrix}}{\mu\mSigma}{\I_{K+1}}$. In the latter case it has been shown \cite{Richmond96c,Richmond96d} that the GLRT is still Kelly's detector \cite{Kelly86} and that its distribution under  the null hypothesis is the same as in the Gaussian case. The assumption here is different since we have a distribution mismatch between $\vx$ and $\X$ which can be direct or the consequence of a particular covariance mismatch.  The aim of the present section is to derive statistical representations of $(\beta,\ttilde)$ under this framework in order to figure out how they deviate from the Gaussian case. 

Let us start with
\begin{align}
	s_{1} &=\vx^{H}\S^{-1}\vx \nonumber \\
	& \dist \mu^{-1} \vx^{H}	\mSigma^{-\onehalf} \F \mSigma^{-\onehalf}  \vx \dist \mu^{-1} \vxtilde^{H} \F \vxtilde
\end{align} 
where $\vxtilde= \Q^{H} \mSigma^{-\onehalf}\vx \dist \vCN{N}{\alpha(\vv^{H}\mSigma^{-1}\vv)^{1/2}}{\I_{N}}$. Similarly
\begin{align}
	s_{2} &=  \frac{|\vx^{H}\S^{-1}\vv|^{2}}{\vv^{H}\S^{-1}\vv} \dist	\mu^{-1} \frac{|\vxtilde^{H}\F\elast|^{2}}{\elast^{H}\F\elast}
\end{align}
Partitioning $\vxtilde=\begin{bmatrix} \vxtilde_{1} \\ \xtilde_{2} \end{bmatrix}$ and $\F$ as in \eqref{partition_F}, it is readily shown that
\begin{align}
	s_{1} &= s_{2} + \mu^{-1} \vxtilde_{1}^{H} \F_{1.2} \vxtilde_{1} \nonumber \\
	s_{2} &= \mu^{-1} F_{22} \left| \xtilde_{2} + \vxtilde_{1}^{H}\vt_{12} \right|^{2}
\end{align}
where $\F_{1.2} = \F_{11} - \F_{12}F_{22}^{-1}\F_{21}$. Consequently
\begin{equation}
	\beta = (1+	\mu^{-1}\vxtilde_{1}^{H} \F_{1.2} \vxtilde_{1})^{-1}
\end{equation}
From  \ref{app:F}, we have that $\F_{1.2} \dist \CF{N-1}{\nu-1}{K} \dist \W_{1}^{\onehalf} \W_{2}^{-1} \W_{1}^{\onehalf}$ where $\W_{1} \dist \CW{N-1}{\nu-1}{\I_{N-1}}$ and $\W_{2} \dist \CW{N-1}{K}{\I_{N-1}}$. Using well-known results on quadratic forms in Wishart distributions, it comes
\begin{equation}
	\vxtilde_{1}^{H} \F_{1.2} \vxtilde_{1} \dist (\vxtilde_{1}^{H}  \vxtilde_{1}) \frac{\Cchisquare{\nu-1}{0}}{\Cchisquare{K-N+2}{0}}
\end{equation}
and finally
\begin{equation}\label{rep_beta_Student}
	\beta \dist \left[1 + \frac{\Cchisquare{\nu-1}{0}}{\mu} \frac{\Cchisquare{N-1}{0}}{\Cchisquare{K-N+2}{0}} \right]^{-1}
\end{equation}
where we used the fact that $\vxtilde_{1} \dist \vCN{N-1}{\vzero}{\I_{N-1}}$ and hence $\vxtilde_{1}^{H}  \vxtilde_{1} \dist \Cchisquare{N-1}{0}$. The previous equation should be compared to its Gaussian counterpart namely
\begin{equation}\label{rep_beta_Gaussian}
	\betaGaussian \dist \left[1 + \frac{\Cchisquare{N-1}{0}}{\Cchisquare{K-N+2}{0}} \right]^{-1}
\end{equation}
The difference lies in the factor $\mu^{-1}\Cchisquare{\nu-1}{0}$. The latter is gamma distributed with mean $\mu^{-1}(\nu-1)$ and variance $\mu^{-2}(\nu-1)$. Note that if one imposes $K^{-1}\E{\X\X^{H}}=(\nu-N)^{-1}\mu \mSigma = \mSigma$ then the mean and variance become $(\nu-1)/(\nu-N)$ and $(\nu-1)/(\nu-N)^{2}$. Therefore as $\nu \rightarrow \infty$ the distribution of this variable becomes more and more concentrated around $1$ and the two representations are equivalent. However, for small $\nu$ there is a difference which will be quantified below. Another observation is that in the Gaussian case $\beta$ and $\rho$ have the same distribution which is no longer the case with Student distributed samples.

\begin{figure}[h]
	\centering
	\includegraphics[width=11cm]{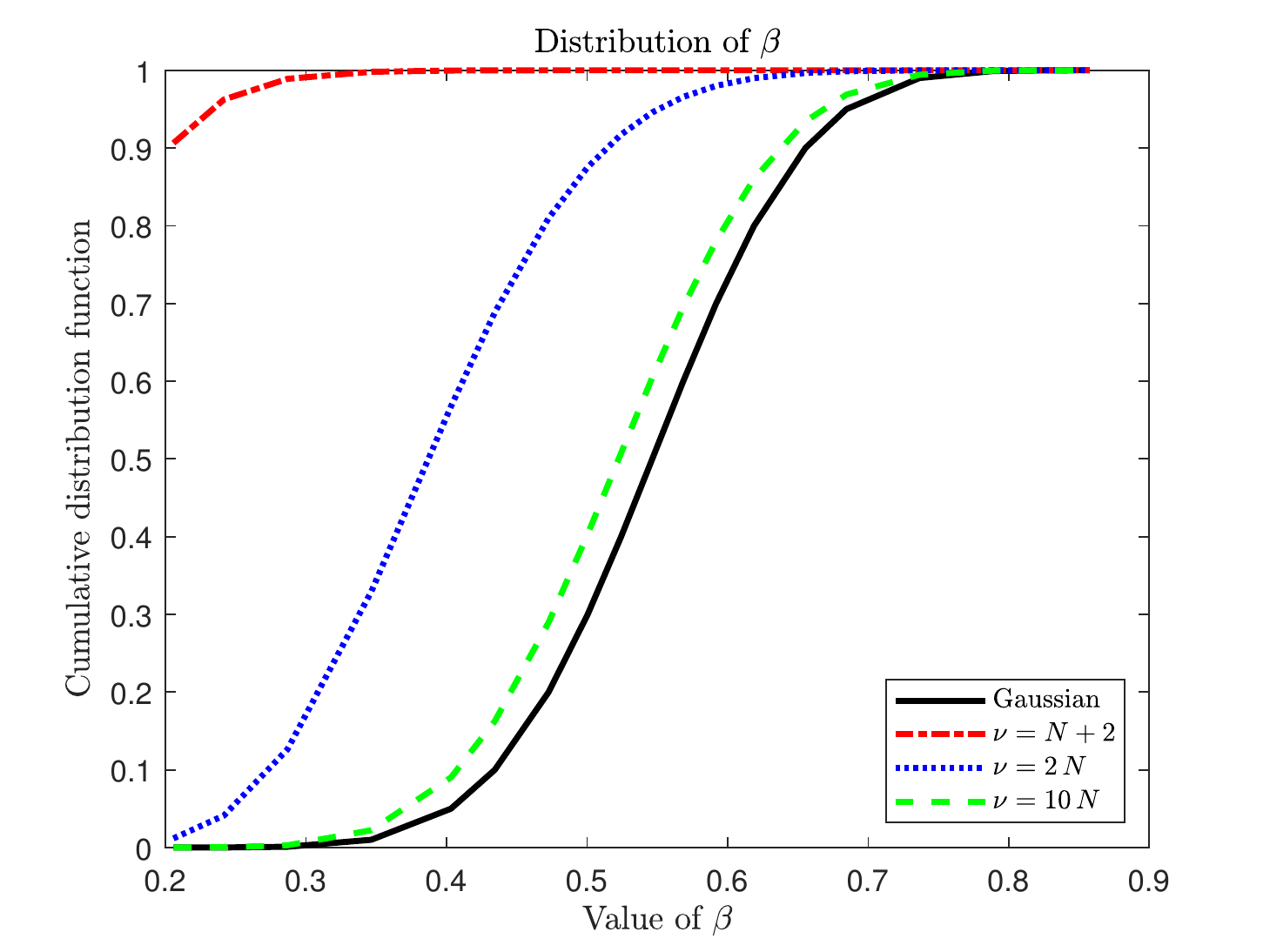}
	\caption{Cumulative distribution function of $\beta$ for various $\nu$. $K=2N$.}
	\label{fig:CDF_beta_K=32_mean_iW=1}
\end{figure}

\begin{figure}[h]
	\centering
	\includegraphics[width=11cm]{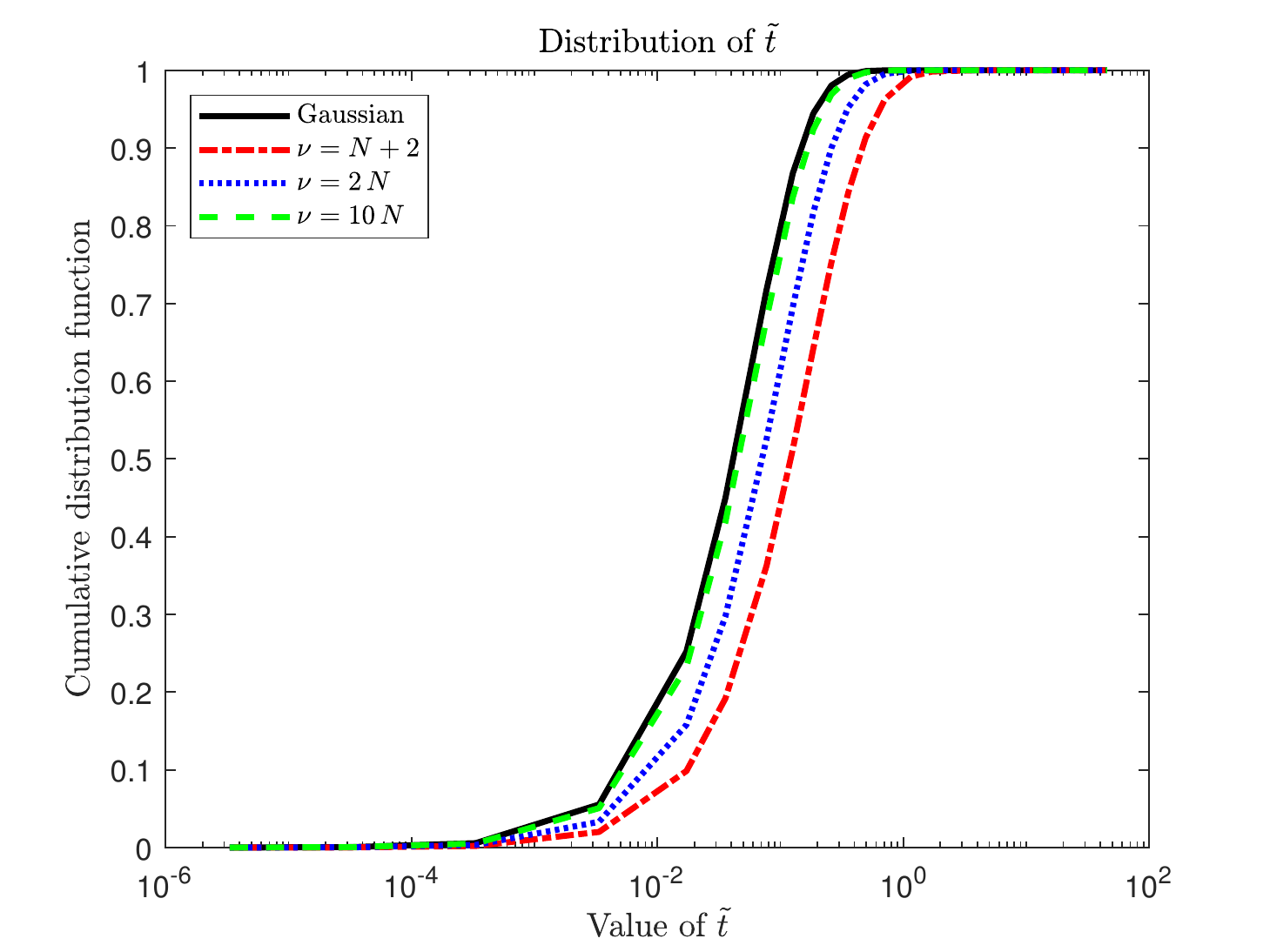}
	\caption{Cumulative distribution function of $\ttilde$ for various $\nu$. $K=2N$.}
	\label{fig:CDF_ttilde_K=32_mean_iW=1}
\end{figure}

\begin{figure}[h]
	\centering
	\includegraphics[width=11cm]{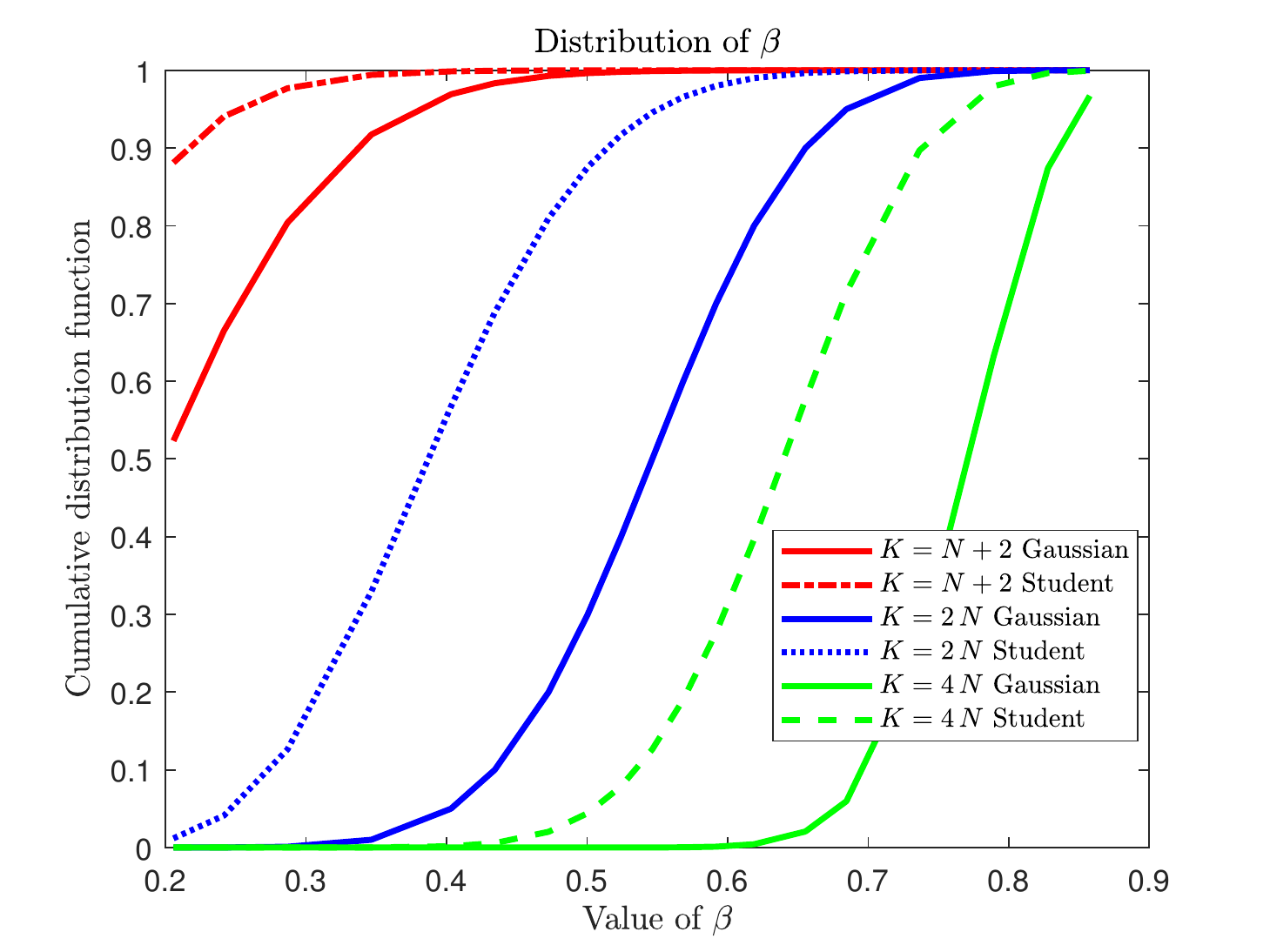}
	\caption{Cumulative distribution function of $\beta$ for various $K$. $\nu=2N$.}
	\label{fig:CDF_beta_nu=32_mean_iW=1}
\end{figure}

\begin{figure}[h]
	\centering
	\includegraphics[width=11cm]{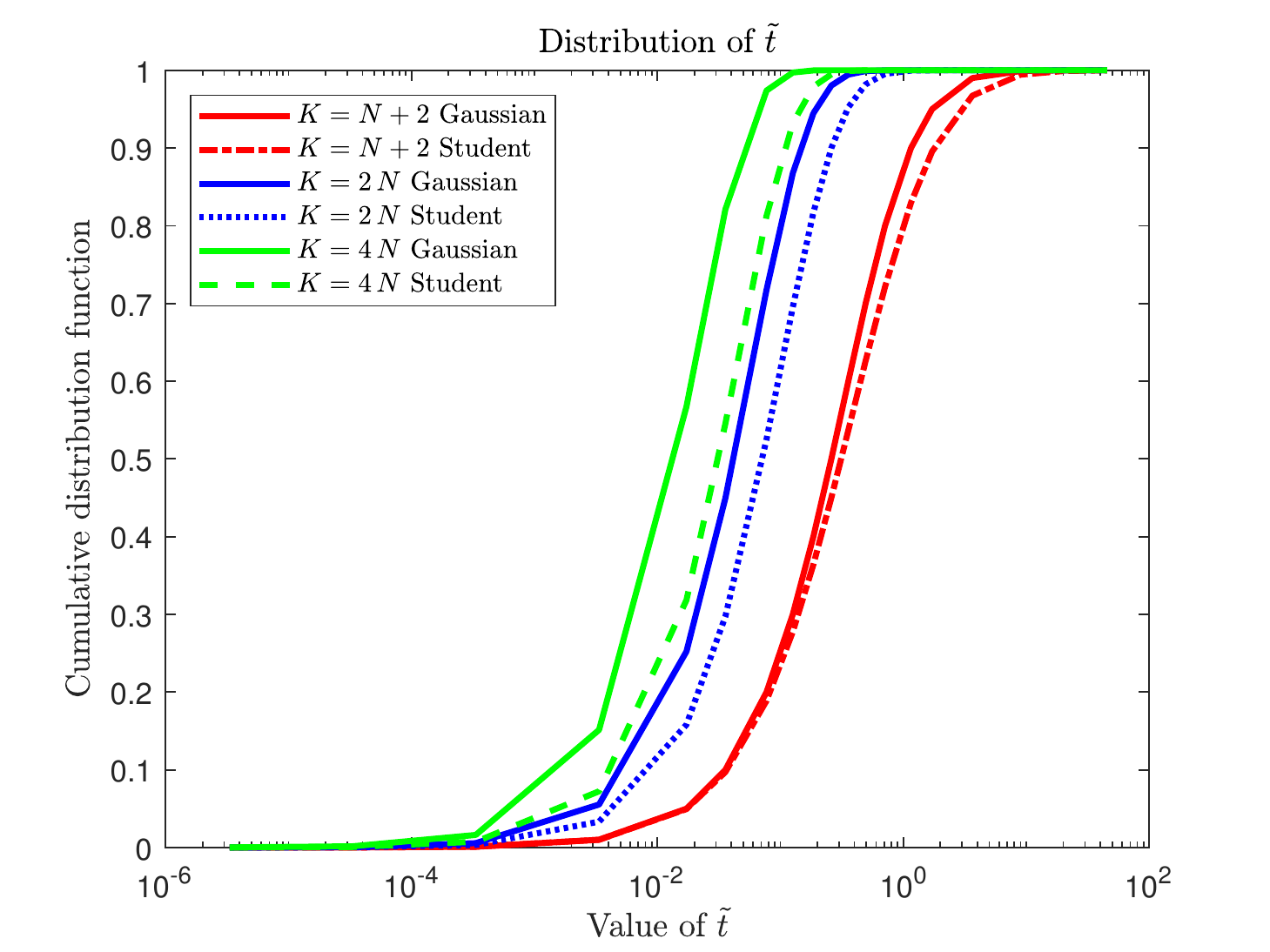}
	\caption{Cumulative distribution function of $\ttilde$ for various $K$. $\nu=2N$.}
	\label{fig:CDF_ttilde_nu=32_mean_iW=1}
\end{figure}

Let us now turn to $\ttilde$ which is the test statistic of Kelly's GLRT and can be written as
\begin{align}
	\ttilde &= \frac{s_{2}}{1+s_{1}-s_{2}} \dist \frac{\mu^{-1} F_{22} \left| \xtilde_{2} + \vxtilde_{1}^{H}\vt_{12} \right|^{2}}{1+\mu^{-1} \vxtilde_{1}^{H} \F_{1.2} \vxtilde_{1}}	
\end{align}
From the representation of $\vt_{12}$ in \eqref{rep_t12}, we have that
\begin{equation}
	\xtilde_{2} + \vxtilde_{1}^{H}\vt_{12} = \xtilde_{2} + (1+F_{22}^{-1})^{1/2} \frac{\vxtilde_{1}^{H}\vn_{12}}{\sqrt{\gamma_{12}}} 
\end{equation}
which implies, since $\xtilde_{2} \dist \CN{\alpha(\vv^{H}\mSigma^{-1}\vv)^{1/2}}{1}$ that
\begin{equation}
	\xtilde_{2} + \vxtilde_{1}^{H}\vt_{12} | \vxtilde_{1},F_{22},\gamma_{12}  \dist \CN{\alpha(\vv^{H}\mSigma^{-1}\vv)^{1/2}}{1+(1+F_{22}^{-1}) \frac{\vxtilde_{1}^{H}\vxtilde_{1}}{\gamma_{12}}}
\end{equation}
Consequently
\begin{align}\label{rep_ttilde_Student}
	\ttilde	\dist \frac{\mu^{-1} F_{22}\left[1+(1+F_{22}^{-1}) \frac{\vxtilde_{1}^{H}\vxtilde_{1}}{\gamma_{12}}\right]}{1 + \frac{\Cchisquare{\nu-1}{0}}{\mu} \frac{\vxtilde_{1}^{H}\vxtilde_{1}}{\Cchisquare{K-N+2}{0}}} \Cchisquare{1}{\delta}
\end{align}
where $\delta =|\alpha|^{2}(\vv^{H}\mSigma^{-1}\vv)/\left[1+(1+F_{22}^{-1}) \frac{\vxtilde_{1}^{H}\vxtilde_{1}}{\gamma_{12}}\right]$, the distributions of $F_{22}$, $\gamma_{12}$ are given in \eqref{rep_F22} and $\vxtilde_{1}^{H}  \vxtilde_{1} \dist \Cchisquare{N-1}{0}$. Equation \eqref{rep_ttilde_Student} provides the statistical representation of $\ttilde$ as a function of independent chi-square distributed random variables. It should be compared with the Gaussian expression
\begin{equation}\label{rep_ttilde_Gaussian}
	\ttildeGaussian | \betaGaussian\dist 	\frac{\Cchisquare{1}{\betaGaussian |\alpha|^{2}(\vv^{H}\mSigma^{-1}\vv)}}{\Cchisquare{K-N+1}{0}}
\end{equation}
We now evaluate how the distributions of $\beta$ and $\ttilde$ in the Student case depart from their distributions with Gaussian distributed training samples. As before we have $N=16$ and $\mu=\nu-N$. Figures \ref{fig:CDF_beta_K=32_mean_iW=1}-\ref{fig:CDF_ttilde_K=32_mean_iW=1}  display the c.d.f. of  $\beta$ and $\ttilde$  for $K=2N$. Similarly to what was observed for the SNR loss we see that the impact is significant, especially for $\beta$. This suggests that using $\beta$ in the test statistic may lead to significant performance degradation. We also notice that contrary to the Gaussian case $\loss$ and $\beta$ do not have the same distribution in the Student case. Furthermore, similarly to what was observed for $\loss$, the difference between Gaussian and Student distributions is all the more important that $K$ is large, see Figures \ref{fig:CDF_beta_nu=32_mean_iW=1}-\ref{fig:CDF_ttilde_nu=32_mean_iW=1}.

\begin{figure}[h]
	\centering
	\includegraphics[width=11cm]{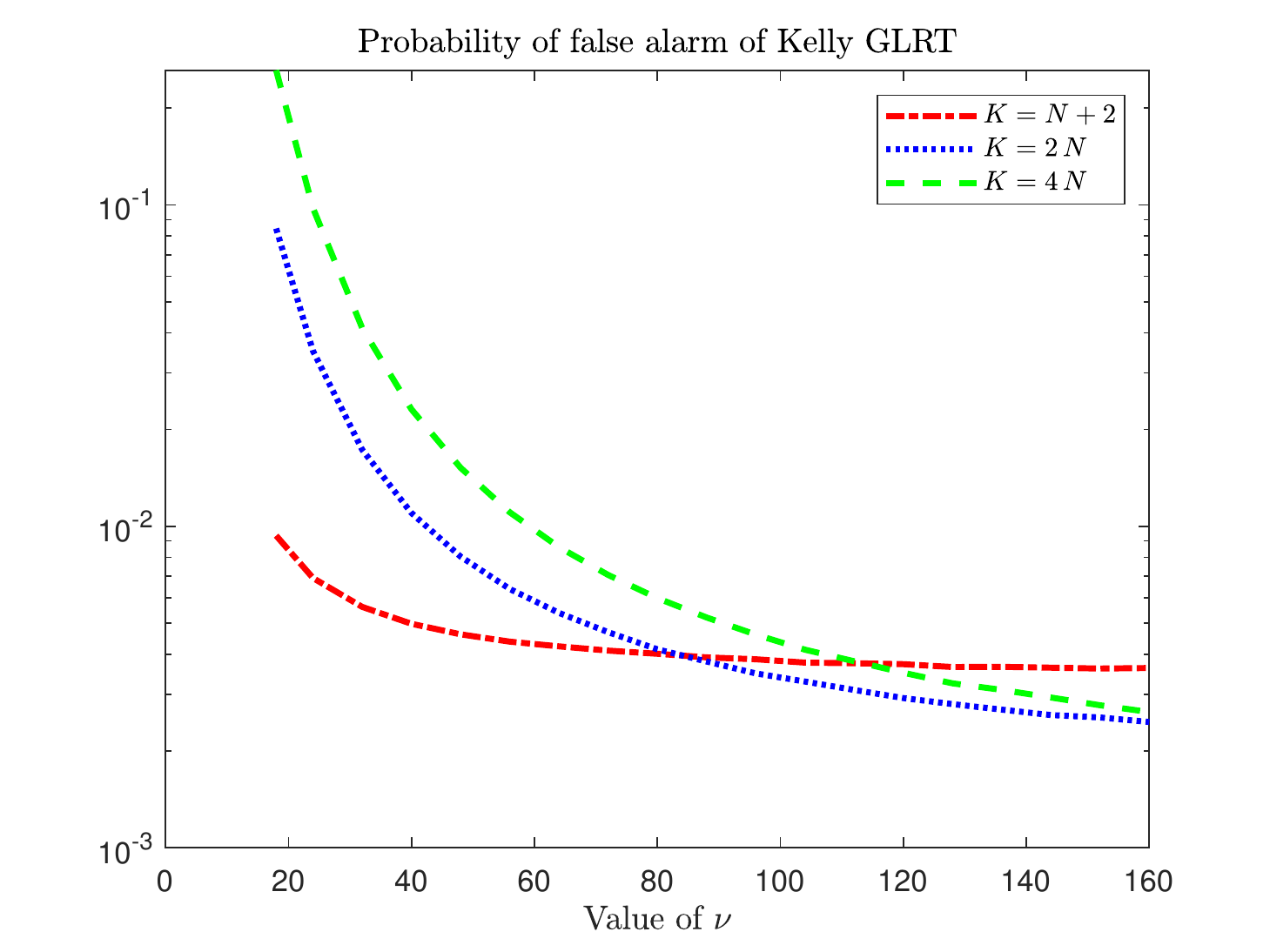}
	\caption{Probability of false alarm of $\ttilde$ versus $\nu$ for various values of $K$.}
	\label{fig:Pfa_Student_mean_iW=1}
\end{figure}

Finally we investigate the influence of Student distributed training samples on the probability of false alarm of $\ttilde$, i.e., Kelly's Gaussian GLRT. The threshold is set so that $\Pfa=10^{-3}$ in the Gaussian case. Figure \ref{fig:Pfa_Student_mean_iW=1} shows the actual $\Pfa$ when $t$ distributed training samples are used. One can observe two things. First, $\Pfa$ is increased and the increase is more pronounced as $K$ grows. Second, one can see that even for large $\nu$ we do not recover the Gaussian $\Pfa$ due to the distribution mismatch between the data under test $\vx$ and the training samples $\X$.
\section{Conclusions}
In this paper we were interested in what happens to statistics commonly used in adaptive multichannel processing when the training samples used to infer noise are no longer Gaussian distributed but $t$ distributed. Statistical representations of the SNR loss and of some statistics used for adaptive detection were derived, based on properties of partitioned matrix-variate $F$ distributions. The expressions derived are given in terms of independent chi-square distributed random variables. They enable one to quickly evaluate the impact of this type of distribution mismatch, which was illustrated numerically.  
\appendix
\section{Properties of partitioned complex matrix-variate F distributed matrices\label{app:F}}
In this appendix we derive some properties of  partitioned complex matrix-variate $F$ distributed matrices. Most of these properties were derived in the real case in \cite{Tan69}. We extend them to the complex case and provide new additional results concerning marginalization of the distribution of $\vt_{12}$, see below. Let $\F \dist \CF{p}{q}{n}$ and let us partition it as
\begin{equation}
	\overset{\begin{matrix} r & s \end{matrix}}{ \begin{pmatrix} \F_{11} & \F_{12} \\ \F_{21} & \F_{22} \end{pmatrix}}	\begin{matrix} r \\ s \end{matrix}
\end{equation}
The p.d.f. of $\F$ is given by $p(\F) \propto \det{\F}^{q-p} \det{\I_{p}+\F}^{-(q+n)}$. Now we have $\det{\F} = \det{\F_{1.2}} \det{\F_{22}}$ where $\F_{1.2} = \F_{11} - \F_{12}\F_{22}^{-1}\F_{21}$.
Moreover
\begin{align}
	\left[\I_{p}+\F\right]_{1.2} &= \I_{r} + \F_{11} - \F_{12}(\I_{s}+\F_{22})^{-1}\F_{21} \nonumber \\
	&= \I_{r} + \F_{1.2} + \F_{12}\F_{22}^{-1}\F_{21} - \F_{12}(\I_{s}+\F_{22})^{-1}\F_{21} \nonumber \\
	&= \I_{r} + \F_{1.2} + \F_{12} \left[\F_{22}^{-1} - (\I_{s}+\F_{22})^{-1}\right]\F_{21} \nonumber \\
	&= \I_{r} + \F_{1.2} + \F_{12} \F_{22}^{-1} (\I_{s}+\F_{22})^{-1}\F_{21} \nonumber \\
	&= \I_{r} + \F_{1.2} + \T_{12}  (\I_{s}+\F_{22})^{-1}\F_{22}\T_{12}^{H}
\end{align}
with $\T_{12}=\F_{12} \F_{22}^{-1}$. It ensues that
\begin{align}
	\det{\left[\I_{p}+\F\right]_{1.2}}	&=\det{\I_{r} + \F_{1.2}} \nonumber \\  
	&\times \det{\I_{r}+(\I_{r} + \F_{1.2})^{-1} \T_{12}  (\I_{s}+\F_{22})^{-1}\F_{22}\T_{12}^{H}}
\end{align}
Since the Jacobian $J(\F \rightarrow \F_{1.2},\T_{12},\F_{22}) = \det{\F_{22}}^{r}$ \cite{Khatri65}, we can write the joint density of $(\F_{1.2},\T_{12},\F_{22})$ as
\begin{align}
	&p(\F_{1.2},\T_{12},\F_{22}) \propto  \det{\F_{1.2}}^{q-r-s} \det{\I_{r} + \F_{1.2}}^{-(q+n-s)} \nonumber \\
	&\times \det{\F_{22}}^{q-s} \det{\I_{s} + \F_{22}}^{-(q+n-r)} \nonumber \\
	&\times \det{\I_{r} + \F_{1.2}}^{-s} \det{\F_{22}^{-1} (\I_{s}+\F_{22})}^{-r} \nonumber \\
	&\det{\I_{r}+(\I_{r} + \F_{1.2})^{-1} \T_{12}  (\I_{s}+\F_{22})^{-1}\F_{22}\T_{12}^{H}}^{-(q+n)}
\end{align}
Therefore  $\F_{1.2}$ and $\F_{22}$ are independent and
\begin{align}\label{pdf_F1.2_F22}
	\F_{1.2}  \dist \CF{r}{q-s}{n}; \quad \F_{22}  \dist \CF{s}{q}{n-r}
\end{align}
\begin{equation}
	\T_{12} | 	\F_{1.2},\F_{22} \dist \mCT{r,s}{n+q-p+1}{\Mzero}{\I_{r} + \F_{1.2}}{\F_{22}^{-1} (\I_{s}+\F_{22})}
\end{equation}
These results extend those of \cite{Tan69} to the complex case. Next, we  marginalize $\T_{12}$ in order to obtain the distribution of $\T_{12}|\F_{22}$. To do so, note that
\begin{align}
	&\det{\I_{r}+(\I_{r} + \F_{1.2})^{-1} \T_{12}  (\I_{s}+\F_{22})^{-1}\F_{22}\T_{12}^{H}} \nonumber \\
	&= \det{\I_{r} + \F_{1.2}}^{-1} \det{\mDelta+\F_{1.2}}
\end{align}
where $\mDelta =\I_{r} + \T_{12}  (\I_{s}+\F_{22})^{-1}\F_{22}\T_{12}^{H}$.  It follows that
\begin{align}
	&p(\T_{12}|\F_{22} ) = \int p(\T_{12}|\F_{22},\F_{1.2}) p(\F_{1.2}) \, \mathrm{d}\F_{1.2} \nonumber \\
	&= \det{\F_{22}^{-1} (\I_{s}+\F_{22})}^{-r} \int \det{\I_{r} + \F_{1.2}}^{n+q-s} \det{\mDelta+\F_{1.2}}^{-(q+n)} p(\F_{1.2}) \, \mathrm{d}\F_{1.2} \nonumber \\
	&= \det{\F_{22}^{-1} (\I_{s}+\F_{22})}^{-r} \int \det{\F_{1.2}}^{q-r-s} \det{\mDelta+\F_{1.2}}^{-(q+n)} \, \mathrm{d}\F_{1.2} \nonumber \\
	&= \det{\F_{22}^{-1} (\I_{s}+\F_{22})}^{-r} \det{\I_{r} + \T_{12}  (\I_{s}+\F_{22})^{-1}\F_{22}\T_{12}^{H}}^{-(n+s)}
\end{align}
This proves that
\begin{equation}\label{pdf_T12|F22}
	\T_{12}|\F_{22} \dist \mCT{r,s}{n-r+1}{\Mzero}{\I_{r}}{\F_{22}^{-1} (\I_{s}+\F_{22})}	
\end{equation}
When $s=1$, \eqref{pdf_T12|F22} reduces to
\begin{equation}\label{pdf_t12|F22}
	\vt_{12}|F_{22}	\dist \vCT{p-1}{n-p+2}{\vzero}{(1+F_{22}^{-1})\I_{p-1}}
\end{equation}
which means that $\vt_{12}$ can be modelled as
\begin{equation}\label{rep_t12|F22}
	\vt_{12} = 	(1+F_{22}^{-1})^{1/2} \frac{\vn_{12}}{\sqrt{\gamma_{12}}}
\end{equation}
with $\vn_{12} \dist \vCN{p-1}{\vzero}{\I_{p-1}}$ and $\gamma_{12} \dist \Cchisquare{n-p+2}{0}$. The distribution of $\vt_{12}$ can be evaluated by marginalizing $p(\vt_{12}|F_{22})$, which gives
\begin{align}
	&p(\vt_{12}) = \int_{0}^{\infty} p(\vt_{12}|F_{22}) p(F_{22}) \, \mathrm{d}F_{22}	 \nonumber \\
	&= \frac{\Gamma(n+1)}{\pi^{p-1}\Gamma(n-p+2)} \nonumber \\
	& \int_{0}^{\infty}  \frac{F_{22}^{p-1}}{(1+F_{22})^{p-1}} [1+F_{22}(1+F_{22})^{-1}\vt_{12}^{H}\vt_{12}]^{-(n+1)} p(F_{22}) \, \mathrm{d}F_{22} \nonumber \\
\end{align}
From \eqref{pdf_F1.2_F22}, $F_{22}$ has a complex scalar $\tilde{F}(q,n-p+1)$ distribution so that
\begin{equation}\label{pdf_F22}
	p(F_{22})	= \frac{F_{22}^{q-1} (1+F_{22})^{-(q+n-p+1)}}{\Beta{q}{n-p+1}} 
\end{equation}
where $\Beta{a}{b}=\Gamma(a)\Gamma(b)/\Gamma(a+b)$. It follows that
\begin{align}\label{pdf_t12}
	p(\vt_{12}) &= \int_{0}^{\infty} \frac{C \, F_{22}^{p+q-2}}{(1+F_{22})^{q+n}} [1+F_{22}(1+F_{22})^{-1}\vt_{12}^{H}\vt_{12}]^{-(n+1)} \, \mathrm{d}F_{22} \nonumber \\
	&=  \frac{C}{\Beta{p+q-1}{n-p+1}} {}_{2}F_{1}(n+1,p+q-1;n+q;-\vt_{12}^{H}\vt_{12})
\end{align}
with $C=\frac{\Gamma(n+1)}{\pi^{p-1}\Gamma(n-p+2)} \frac{1}{B(q,n-p+1)}$ and where we used \cite{Gradshteyn07} to obtain the last line. The unconditional distribution is seen to depend only on $\vt_{12}^{H}\vt_{12}$. Finally note that in the purpose of analyzing the SNR loss the conditional distribution $p(\vt_{12}|F_{22})$ is the most convenient and is actually used.
\section{Distribution and average value of SNR loss in the Student case \label{app:mean_SNRloss_Student}}
In this appendix we derive the p.d.f as well as the mean value of the SNR loss. Let 	$F_{1} \dist \Cchisquare{K-N+1}{0} / \Cchisquare{\nu}{0}$ and $F_{2} \dist \Cchisquare{N-1}{0} / \Cchisquare{K-N+2}{0}$ and let $\loss \dist [1+(1+F_{1})F_{2}]^{-1}$. Let us first evaluate the distribution of $\loss|F_{1}$. The p.d.f of $F_{2}$ is given by
\begin{equation}\label{p(F2)}
	p_{F_{2}}(f_{2}) =\frac{1}{\Beta{N-1}{K-N+2}} \frac{f_{2}^{N-2}}{(1+f_{2})^{K+1}}
\end{equation}
Making the change of variables $\loss = [1+(1+F_{1})F_{2}]^{-1} \Leftrightarrow F_{2}=(1+F_{1})^{-1}(\loss^{-1}-1)$ whose Jacobian is $J(F_{2} \rightarrow \loss|F_{1})= (1+F_{1})^{-1}\loss^{-2}$, it follows that
\begin{equation}\label{p(rho|F1)}
	p(\rho|F_{1}) = 	\frac{(1+F_{1})^{K-N+2}}{\Beta{N-1}{K-N+2}} \frac{\loss^{K-N+1} (1-\loss)^{N-2}}{(1+\loss F_{1})^{K+1}}
\end{equation}
Setting $F_{1}=0$, one recovers the usual beta distribution of the SNR loss in the Gaussian case. Marginalizing with respect to the p.d.f. of $F_{1}$ we obtain
\begin{align}\label{p(rho)}
	p(\loss) &= \int_{0}^{\infty}	p(\loss|f_{1}) p_{F_1}(f_{1}) \, \mathrm{d}f_{1} \nonumber \\
	&= \frac{\loss^{K-N+1} (1-\loss)^{N-2}}{\Beta{N-1}{K-N+2}\Beta{K-N+1}{\nu}} \int_{0}^{\infty} \frac{f_{1}^{K-N}(1+f_{1})^{-(\nu-1)}}{(1+\loss f_{1})^{K+1}} \, \mathrm{d}f_{1} \nonumber \\
	&=  \frac{\Beta{K-N+1}{\nu+N-1}}{\Beta{N-1}{K-N+2}\Beta{K-N+1}{\nu}} \loss^{K-N+1} (1-\loss)^{N-2} {\;}_{2}F_{1}(K+1,K-N+1;\nu+K;1-\loss)
\end{align}
where we made use of \cite{Gradshteyn07} to obtain the last equality. The previous equation allows to calculate the average value of the SNR loss. Let us start with the conditional mean of $\loss$:
\begin{align}
	\E{\loss | F_{1}} &=  \int_{0}^{1} \loss \; p(\rho|F_{1}) \, \mathrm{d}\loss \nonumber \\
	&= \frac{(1+F_{1})^{K-N+2}}{\Beta{N-1}{K-N+2}}  \int_{0}^{1}\frac{\loss^{K-N+2} (1-\loss)^{N-2}}{(1+\loss F_{1})^{K+1}} \, \mathrm{d}\loss \nonumber \\
	&=  \frac{\Beta{N-1}{K-N+3}}{\Beta{N-1}{K-N+2}} (1+F_{1})^{K-N+2} {\;}_{2}F_{1}(K+1,K-N+3;K+2;-F_{1}) \nonumber \\
	&=  \frac{K-N+2}{K+1} (1+F_{1})^{K-N+2} {\;}_{2}F_{1}(K+1,K-N+3;K+2;-F_{1}) \nonumber \\
	&=  \frac{K-N+2}{K+1} (1+F_{1})^{-1} {\;}_{2}F_{1}\left(1,K-N+3;K+2;\frac{F_{1}}{1+F_{1}}\right) 
\end{align}
where the two last lines are obtained from equivalent expressions of the hypergeometric function \cite{Gradshteyn07}. If we set $F_{1}=0$ in the previous equation we recover the Gaussian case for which $\E{\lossGaussian}=\frac{K-N+2}{K+1}$. Next we need to integrate with respect to the density of $F_{1}$:
\begin{align}
	\E{\loss} &= \int_{0}^{\infty} \E{\loss | f_{1}} p_{F_{1}} (f_{1}) \, \mathrm{d}f_{1}\nonumber \\
	& = \int_{0}^{\infty}	\frac{\E{\loss | f_{1}} }{\Beta{K-N+1}{\nu}}  \frac{f_{1}^{K-N}}{(1+f_{1})^{\nu+K-N+1}} \, \mathrm{d}f_{1} \nonumber \\
	&= \frac{(K-N+2)}{(K+1)\Beta{K-N+1}{\nu}}  \int_{0}^{\infty} \frac{f_{1}^{K-N}}{(1+f_{1})^{\nu+K-N+2}} {\;}_{2}F_{1}\left(1,K-N+3;K+2;\frac{f_{1}}{1+f_{1}}\right) \, \mathrm{d}f_{1} 
\end{align}
Making the change of variables $x=f_{1}/(1+f_{1})$ the integral above can be written as
\begin{align}
	I &= \int_{0}^{\infty} x^{K-N} (1-x)^{\nu} {\;}_{2}F_{1}(1,K-N+3;K+2;x) \, \mathrm{d}x \nonumber \\
	&= \Beta{K-N+1}{\nu+1}   {\;}_{3}F_{2}(1,K-N+3,K-N+1;K+2,\nu+K-N+2;1)
\end{align}
which finally results in
\begin{align}
	&\E{\loss } = \frac{K-N+2}{K+1} \frac{\Beta{K-N+1}{\nu+1}}{\Beta{K-N+1}{\nu}} \nonumber \\
	&\times {}_{3}F_{2}(1,K-N+3,K-N+1;K+2,\nu+K-N+2;1) \nonumber \\
	&= \frac{\nu(K-N+2)}{(\nu+K-N+1)(K+1)} \nonumber \\
	&\times {}_{3}F_{2}(1,K-N+3,K-N+1;K+2,\nu+K-N+2;1)
\end{align}

\end{document}